\documentclass[12pt]{article}
\usepackage[british]{babel}
\usepackage{amsmath,amsthm,amssymb,bm,amsfonts,enumerate,enumitem}
\setenumerate[1]{itemsep=2pt,partopsep=0pt,parsep=\parskip,topsep=5pt}
\setitemize[1]{itemsep=2pt,partopsep=0pt,parsep=\parskip,topsep=5pt}
\usepackage{graphicx, subfigure,geometry,indentfirst,setspace,authblk}
\usepackage[font=scriptsize]{caption}
\usepackage{xcolor}
\usepackage{float}
\usepackage{hyperref}
\makeatletter
  \newcommand\figcaption{\def\@captype{figure}\caption}
  \newcommand\tabcaption{\def\@captype{table}\caption}
  \newcommand{\RomanNum}[1]{\uppercase\expandafter{\romannumeral#1}}

\makeatother
\newtheoremstyle{special}
    {\topsep}
    {\topsep}
    {\itshape}
    {}
    {\bfseries}
    {}
    {0.5em}
    {{\thmname{#1}\thmnumber{ #2$^{\bm*}\!$.}\thmnote{\ \textmd{(#3)}}}}

\newtheorem{thm}{Theorem}[section]
\newtheorem{cor}[thm]{Corollary}
\newtheorem{conj}[thm]{Conjecture}
\newtheorem{lem}[thm]{Lemma}

\newtheorem{obs}[thm]{Observation}

\newcounter{casecounter}

\numberwithin{equation}{section}

\theoremstyle{special}

\hypersetup{hidelinks}
\begin{document}
\title{Quasi-Bisections of Cubic Graphs with Tree and Circuit Decompositions}
\author[1]{Xixiong Tang}
\author[1]{Rong-Xia Hao\footnote {Corresponding author and email: rxhao@bjtu.edu.cn}}
\author[1]{Yilun Luo}
\author[2]{Cun-Quan Zhang}
\renewcommand{\Affilfont}{\small}
\affil[1]{School of Mathematics and Statistics, Beijing Jiaotong University, Beijing 100044, China}
\affil[2]{Department of Mathematics, West Virginia University, Morgantown, WV 26506, USA.}
\date{}
\maketitle
\begin{abstract}
A \(2\)-quasi-bisection of a graph is a vertex \(2\)-coloring in which the color classes differ in size by at most two and every monochromatic component has at most two vertices. Esperet, Mazzuoccolo, and Tarsi in [J. Graph Theory, 86(2017), 149--158] conjectured that every cubic graph admits such a coloring. DeVos and Nurse proved the conjecture for cubic graphs whose edge sets can be partitioned into a tree and a circuit. This result is generalized in this paper for a larger family of cubic graphs. The conjecture is verified for cubic graphs whose edge set can be partitioned into two subsets \(E(T)\) and \(E(J)\), where \(T\) is a tree of \(G\) and \(J\) is an even subgraph containing at most two odd circuit components.
\par\vspace{0.5em}
\noindent\textbf{Keywords}: \(2\)-quasi-bisection; cubic tree; color switching; Ban--Linial conjecture.
\par\noindent\textbf{2020 Mathematics Subject Classification}: 05C15, 05C70.
\end{abstract}
\section{Introduction}

A \emph{\(2\)-bisection} is a vertex \(2\)-coloring in which every monochromatic component has at most two vertices and the color classes have equal size.

\begin{conj}[{\cite{ref3}}]\label{conj:ban-linial}
\emergencystretch=1em
Every bridgeless cubic graph admits a \(2\)-bisection, except for the Petersen graph.
\end{conj}

Ban and Linial \cite{ref3} proved the conjecture for \(3\)-edge-colorable cubic graphs. Abreu et al. proved it for bridgeless claw-free cubic graphs \cite{ref1}. They also proved it for \emph{cycle permutation graphs}, with the Petersen graph as the sole exception \cite{ref2}. The connection with flows gives another sufficient condition: every cubic graph with circular flow number less than five admits a \(2\)-bisection \cite{ref5}. Thus graphs with larger circular flow number are of particular interest. Zerafa \cite{ref8} proved the existence of \(2\)-bisections for treelike snarks, a family of graphs with circular flow number at least five. Li et al. \cite{ref6} proved that several infinite families of Halin snarks, which together include all treelike snarks, admit \(2\)-bisections.

The conclusion of Conjecture~\ref{conj:ban-linial} fails without the bridgelessness assumption. Esperet et al. \cite{ref5} constructed an infinite family of cubic graphs with bridges that admit no \(2\)-bisection. This motivates allowing a small difference between the sizes of the color classes.

Following Esperet et al. \cite{ref5}, a \emph{\(2\)-quasi-bisection} is a vertex \(2\)-coloring in which every monochromatic component has at most two vertices and the color classes differ in size by at most two. They proposed the following conjecture.

\begin{conj}[{\cite{ref5}}]\label{conj:quasi-bisection}
Every cubic graph admits a \(2\)-quasi-bisection.
\end{conj}

DeVos and Nurse \cite{ref4} proved the following result.

\begin{thm}[{\cite{ref4}}]\label{thm:devos-nurse}
Let \(G\) be a cubic graph. Then \(G\) admits a \(2\)-quasi-bisection if either of the following conditions holds:
\begin{enumerate}[label=(\arabic*)]
\item \(G\) contains a cubic tree \(T\) such that \(G-E(T)\) is bipartite;
\item \(E(G)\) can be partitioned into the edge sets of two subgraphs \(T\) and \(J\), where \(T\) is a tree and \(J\) is a circuit.
\end{enumerate}
\end{thm}

We extend the tree-and-circuit case to the following decompositions.

\begin{thm}\label{result:1.1}

Let \(G\) be a cubic graph. Suppose that \(E(G)\) can be partitioned into the edge sets of two subgraphs \(T\) and \(J\), where \(T\) is a tree and \(J\) is an even subgraph containing at most two odd circuit components. Then \(G\) admits a \(2\)-quasi-bisection.

\end{thm}

The proof proceeds in three steps. We first precolor the circuits of \(J\) and then extend this coloring to all vertices of \(T\). Finally, a \(2\)-quasi-bisection is obtained by applying switching operations to a subset of vertices.

Section 2 introduces the notation and two basic facts about cubic trees. Section 3 establishes the tree extension lemma, and Section 4 develops switching operations along bichromatic paths. Section 5 combines these tools to prove Theorem~\ref{result:1.1}.

\section{Preliminaries}

All graphs in this paper are finite and simple. A graph is \emph{cubic} if every vertex is of degree three. A tree is called a \emph{cubic tree} if every vertex is of degree \(1\) or \(3\). In particular, \(K_2\) is a cubic tree. We denote the leaf set of \(T\) by \(L(T)\). A \emph{circuit} is a connected \(2\)-regular subgraph. A subgraph is \emph{even} if every vertex is of even degree in that subgraph.

\begin{obs}\label{obs:even-complement-tree}
Let \(G\) be a cubic graph and \(T\) a tree contained in \(G\). If \(G-E(T)\) is an even subgraph, then \(T\) is a spanning cubic tree of \(G\).
\end{obs}

Let \(J=G-E(T)\). Since \(G\) is cubic and \(J\) is even, every vertex of \(G\) is incident with at least one edge of \(T\). Hence \(T\) spans \(G\). Moreover, \(d_T(v)=3-d_J(v)\in\{1,3\}\) for every \(v\in V(G)\), so \(T\) is a cubic tree.

For a vertex \(x\) of a graph \(H\), \(N_H(x)\) denotes its set of neighbors. For \(X\subseteq V(G)\), the \emph{closed neighborhood} \(N_G[X]\) consists of the vertices in \(X\) and all vertices adjacent to at least one vertex in \(X\). A \emph{chord} of a path \(P\) in \(G\) is an edge of \(G\) joining two nonconsecutive vertices in \(V(P)\).

We use \(B\) and \(W\) to denote the colors black and white, respectively. The same symbols also denote the corresponding color classes. An edge is \emph{monochromatic} if its endpoints have the same color and \emph{bichromatic} otherwise. A \emph{monochromatic component} is a connected component of either induced subgraph \(G[B]\) or \(G[W]\). The \emph{monochromatic degree} of a vertex is its number of neighbors of the same color.

Let \(T\) be a cubic tree and \(x,y\) two distinct leaves of \(T\). We call \(\{x,y\}\) a \emph{pair of cherries} if \(N_T(x)\cap N_T(y)\ne\varnothing\). The unique vertex \(u\in N_T(x)\cap N_T(y)\) is the \emph{center} of \(\{x,y\}\). We use the following fact.

\begin{obs}[{\cite{ref4}}]\label{result:3.1}

Every cubic tree with at least four vertices contains a pair of cherries.

\end{obs}

\section{Extending precolorings of the leaves}

 The following lemma is proved by a similar method developed in  \cite{ref4}.

For any vertex coloring \(\phi\), including a prescribed leaf coloring, let \(D_\phi=|B|-|W|\), where \(B,W\) are its color classes.

\begin{lem}\label{result:3.2}
Let \(T\) be a cubic tree with leaf set \(L\), and let \(f:L\to\{B,W\}\) be a coloring of its leaves. Then there exists a pair of colorings \(f'_{+1},f'_{-1}:V(T)\to\{B,W\}\) such that
\begin{enumerate}
\item[(1)] each coloring preserves the prescribed color \(f(v)\) of every leaf \(v\in L\),
\item[(2)] under each coloring, the monochromatic edges form a matching, and
\item[(3)] the color-class differences satisfy
\begin{equation}\label{eq:tree-bounds}
\frac{2D_f-2}{3}\le D_{f'_{+1}}\le\frac{2D_f+6}{3},
\qquad
\frac{2D_f-6}{3}\le D_{f'_{-1}}\le\frac{2D_f+2}{3}.
\end{equation}
\end{enumerate}
\end{lem}

\begin{proof}
We prove both extensions by induction on \(|V(T)|\).

\textbf{Base case.} \(T=K_2\). Set \(f'_{+1}=f'_{-1}=f\). The monochromatic edges form a matching. All vertices are leaves, so \(D_{f'_{+1}}=D_{f'_{-1}}=D_f\in\{-2,0,2\}\). Thus both bounds in \eqref{eq:tree-bounds} hold.

Assume that \(|V(T)|\ge4\). Suppose that the statement holds for every smaller cubic tree \(T^*\) with any prescribed coloring of its leaves. The lemma is proved in two cases.

\medskip
\textbf{Case 1.} There is a monochromatic pair of cherries. Let \(x,y\) be the leaves and \(u\) their common neighbor. Delete \(x,y\) to obtain a smaller cubic tree \(T^*\). Color its new leaf \(u\) opposite to \(x,y\), and retain all other prescribed leaf colors. Let \(h\) be this leaf coloring of \(T^*\). The two cases are shown in Figure~\ref{fig:case1-cherry}.

\begin{figure}[H]
\centering
\includegraphics[width=0.94\linewidth]{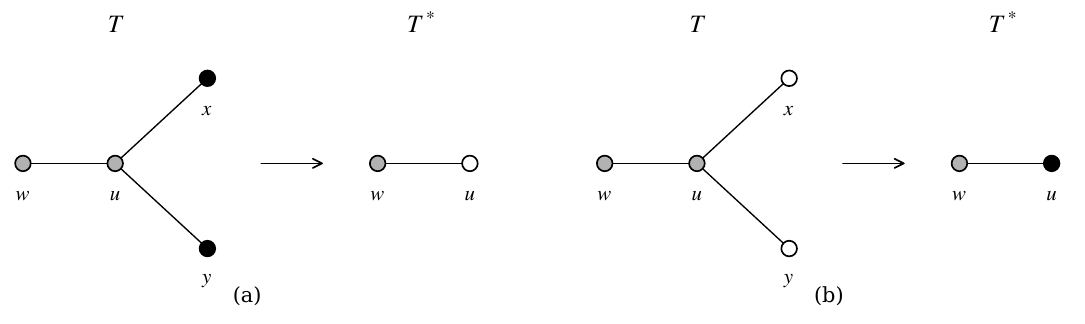}
\caption{Deleting the cherry pair \(\{x,y\}\) in Case 1: (a) black leaves; (b) white leaves. The vertex \(w\) is the remaining neighbor of \(u\). Gray vertices have no color specified in the drawing.}\label{fig:case1-cherry}
\end{figure}

By the induction hypothesis, applied with the two choices \(+1\) and \(-1\), respectively, we obtain extensions \(h'_{+1}\) and \(h'_{-1}\) of \(h\). For each \(\varepsilon\in\{+1,-1\}\), restore \(x,y\) with their prescribed colors and keep \(h'_{\varepsilon}\) on \(V(T^*)\); denote the resulting coloring by \(f'_{\varepsilon}\). Thus \(f'_{+1}\) and \(f'_{-1}\) are the two required extensions. Every prescribed leaf color is preserved. Both restored edges are bichromatic, so the monochromatic edges remain a matching.

We now verify \eqref{eq:tree-bounds}. If \(x,y\) are black, their deletion and the new white leaf give \(D_h=D_f-3\). Restoring these two black vertices gives
\[
D_{f'_{+1}}=D_{h'_{+1}}+2,
\qquad
D_{f'_{-1}}=D_{h'_{-1}}+2.
\]
If \(x,y\) are white, the corresponding equalities are
\[
D_h=D_f+3,
\qquad
D_{f'_{+1}}=D_{h'_{+1}}-2,
\qquad
D_{f'_{-1}}=D_{h'_{-1}}-2.
\]
By the induction hypothesis,
\[
\frac{2D_h-2}{3}\le D_{h'_{+1}}\le\frac{2D_h+6}{3},
\qquad
\frac{2D_h-6}{3}\le D_{h'_{-1}}\le\frac{2D_h+2}{3}.
\]
Substituting the above relations into these two bounds gives respectively
\[
\frac{2D_f-2}{3}\le D_{f'_{+1}}\le\frac{2D_f+6}{3},
\qquad
\frac{2D_f-6}{3}\le D_{f'_{-1}}\le\frac{2D_f+2}{3}.
\]
This proves \eqref{eq:tree-bounds} in Case 1.

\medskip
\textbf{Case 2.} The leaves in every pair of cherries have different colors. No vertex has three leaf neighbors, since two of them would have the same color. Thus every cherry center has exactly two leaf neighbors and one nonleaf neighbor.

To find a suitable reduction, delete the two leaf neighbors of every cherry center and denote the remaining tree by \(S\). Each former center becomes a leaf, and every other surviving vertex keeps its degree. Thus \(S\) is a cubic tree with at least two vertices.

If \(S=K_2\) (see Figure~\ref{fig:cherry-base}), then both vertices are former cherry centers. For otherwise, if one of them were not a former cherry center, then \(T\) would be \(K_2\) or \(K_{1,3}\), both excluded in this case. Color the two centers differently and restore the prescribed leaf colors. The edge of \(S=K_2\) is bichromatic, and each center is incident with exactly one monochromatic leaf edge. The two monochromatic edges form a matching.

\begin{figure}[H]
\centering
\includegraphics[width=0.42\linewidth]{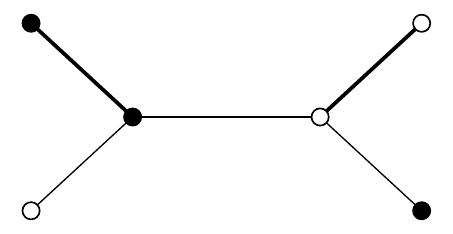}
\captionsetup{justification=centering}\caption{\(S=K_2\).}\label{fig:cherry-base}
\end{figure}

As shown in Figure~\ref{fig:cherry-base}, we have \(D_f=0\), and the resulting coloring has equally sized color classes. This coloring satisfies both bounds in \eqref{eq:tree-bounds}.

We may therefore assume that \(|V(S)|\ge4\). By Observation~\ref{result:3.1}, \(S\) contains a cherry-pair \(\{w_1,w_2\}\). The vertices \(w_1,w_2\) cannot both be original leaves of \(T\): they would form a pair of cherries in \(T\) and would have been deleted. Thus either both are former cherry centers (Subcase 2.1), or exactly one is (Subcase 2.2). The reductions in these two subcases are performed on \(T\), not on \(S\).

\medskip
\textbf{Subcase 2.1.} Both \(w_1,w_2\) are cherry centers of \(T\) (see Figure~\ref{fig:subcase21}). Let \(v\) be their common neighbor. Delete \(w_1,w_2\) and their four leaf neighbors from \(T\) to obtain a smaller cubic tree \(T^*\). The vertex \(v\) is a new leaf of \(T^*\).

We first construct \(f'_{+1}\) (see Figure~\ref{fig:subcase21}(a)). Let \(h\) be a coloring of \(L(T^*)\) in which \(v\) is colored \(W\) and every other leaf retains its color under \(f\). By the induction hypothesis, applied with the two choices \(+1\) and \(-1\), respectively, we obtain extensions \(h'_{+1}\) and \(h'_{-1}\) of \(h\). To construct \(f'_{+1}\), use \(h'_{-1}\): restore the four leaves with their prescribed colors, give \(w_1,w_2\) color \(B\), and keep \(h'_{-1}\) on \(V(T^*)\). Denote this coloring by \(f'_{+1}\).

The leaf colors are preserved. Both edges \(vw_1,vw_2\) are bichromatic, and each restored center is incident with exactly one black monochromatic leaf edge. Hence the monochromatic edges remain a matching.

We next check the first bound in \eqref{eq:tree-bounds}. The four deleted leaves contain two vertices of each color, and the new leaf \(v\) is white. Thus \(D_h=D_f-1\). Restoration adds two black centers and equally many black and white leaves, so \(D_{f'_{+1}}=D_{h'_{-1}}+2\). By the induction hypothesis,
\[
\frac{2D_h-6}{3}\le D_{h'_{-1}}\le\frac{2D_h+2}{3}.
\]
Substituting \(D_h=D_f-1\) and \(D_{h'_{-1}}=D_{f'_{+1}}-2\) into this inequality gives
\[
\frac{2D_f-2}{3}
\le D_{f'_{+1}}
\le \frac{2D_f+6}{3},
\]
which is the first bound in \eqref{eq:tree-bounds}.

To construct \(f'_{-1}\) (see Figure~\ref{fig:subcase21}(b)), let \(h\) now be the leaf coloring of \(T^*\) with \(h(v)=B\) and \(h(x)=f(x)\) for every other leaf \(x\). By the induction hypothesis, take an extension \(h'_{+1}\) of this coloring. Color \(w_1,w_2\) white, restore the four leaves with their prescribed colors, and keep \(h'_{+1}\) on \(V(T^*)\). Here \(D_h=D_f+1\) and \(D_{f'_{-1}}=D_{h'_{+1}}-2\). Substituting these equalities into the induction bound for \(h'_{+1}\) gives
\[
\frac{2D_f-6}{3}
\le D_{f'_{-1}}
\le \frac{2D_f+2}{3}.
\]
This proves \eqref{eq:tree-bounds} in Subcase 2.1.

\begin{figure}[H]
\centering
\includegraphics[width=0.98\linewidth]{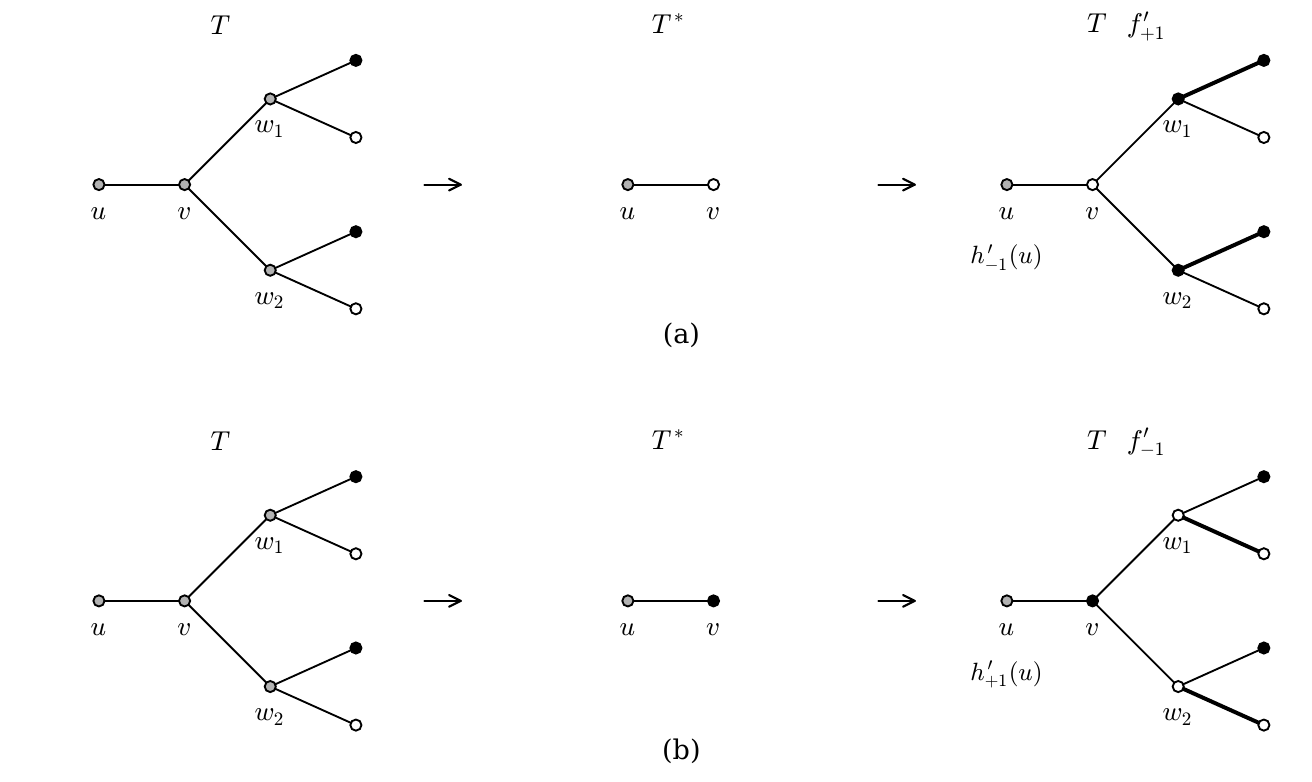}
\caption{The two extensions in Subcase 2.1. Each row shows \(T\) with its prescribed leaf colors, \(T^*\) with its new leaf color, and the extension to \(T\): (a) \(h'_{-1}\) gives \(f'_{+1}\); (b) \(h'_{+1}\) gives \(f'_{-1}\). All vertices of \(T^*\) retain their induction coloring \(h'_{-\varepsilon}\) in the extension \(f'_\varepsilon\). Gray vertices have no color specified in the drawing unless a color label is given. Bold edges are monochromatic.}\label{fig:subcase21}
\end{figure}

\medskip
\textbf{Subcase 2.2.} Exactly one of \(w_1,w_2\) is a cherry center of \(T\). We may assume that \(w_1\) is a leaf of \(T\) and \(w_2\) is a cherry center of \(T\) (see Figure~\ref{fig:subcase22}). Let \(v\) be their common neighbor and \(u\) the third neighbor of \(v\).

Let
\[
T^{**}=\overline{T-\{w_2\}-(N_T(w_2)\cap L(T))}.
\]
This cubic tree is obtained from \(T\) by deleting \(w_2\) and the pair of cherries adjacent to \(w_2\) and suppressing the degree-two vertex \(v\). Let \(h:L(T^{**})\to\{B,W\}\) be a coloring such that \(h(z)=f(z)\) for every \(z\in L(T^{**})\). The leaves of \(T^{**}\) are precisely those of \(T\) except for the two leaf neighbors of \(w_2\). These two leaves have opposite colors, and \(h\) preserves the colors of all remaining leaves. Hence \(D_h=D_f\). By the induction hypothesis, the smaller cubic tree \(T^{**}\) admits a pair of vertex colorings \(h'_{+1},h'_{-1}\) extending \(h\).

We now extend \(h'_{+1}\) and \(h'_{-1}\) to obtain the required colorings of \(T\). For each \(\varepsilon\in\{+1,-1\}\), first define a partial coloring \(g_\varepsilon:\{u,v,w_2\}\to\{B,W\}\) for the path \(uvw_2\) such that \(g_\varepsilon(u)=h'_\varepsilon(u)\) and both edges \(uv,vw_2\) are bichromatic. Then, for each vertex \(z\in V(T)\), define
\begin{equation}\label{eq:subcase22-coloring}
f'_\varepsilon(z)=
\begin{cases}
f(z), & \text{if }z\in L(T);\\
h'_\varepsilon(z), & \text{if }z\in V(T^{**});\\
g_\varepsilon(z), & \text{if }z\in\{v,w_2\}.
\end{cases}
\end{equation}

\begin{figure}[H]
\centering
\includegraphics[width=0.98\linewidth]{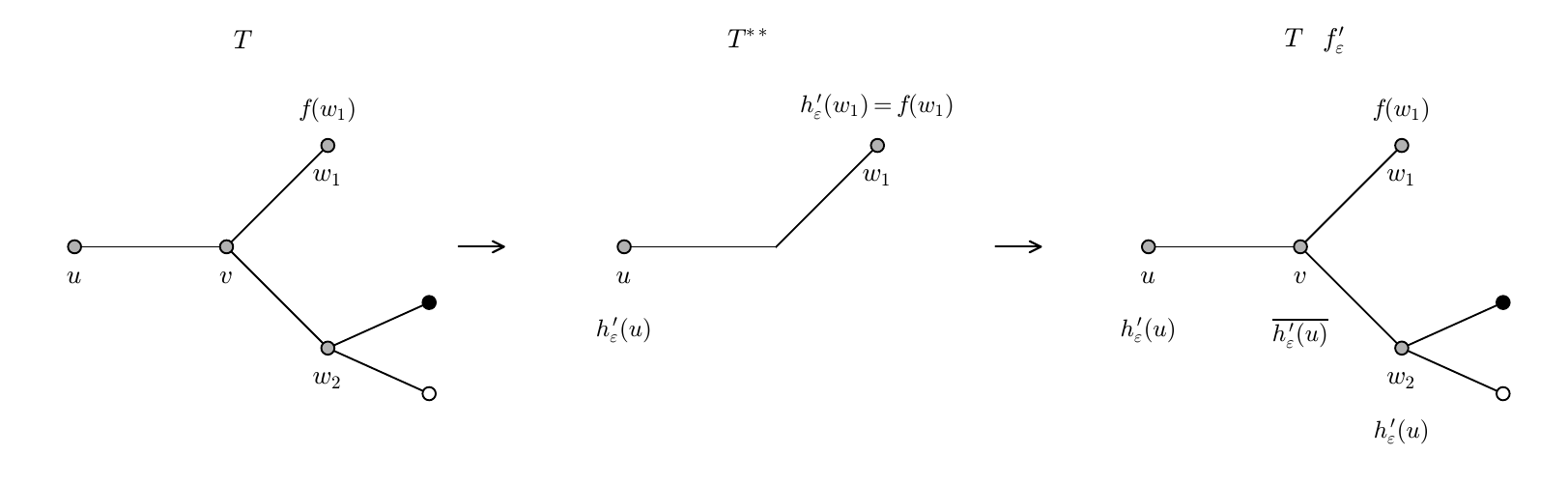}
\caption{Subcase 2.2. The middle drawing shows \(T^{**}\) under \(h'_\varepsilon\), with \(h'_\varepsilon(w_1)=f(w_1)\). The bent line in \(T^{**}\) represents the single edge \(uw_1\); its bend is not a vertex. Gray vertices have no color specified in the drawing unless a color label is given. Both leaves of the cherry pair at \(w_2\) retain their prescribed opposite colors. The bar over \(h'_\varepsilon(u)\) denotes the opposite color.}\label{fig:subcase22}
\end{figure}

By the definition of \(f'_\varepsilon\) (see Equation~\eqref{eq:subcase22-coloring}), each \(f'_\varepsilon\) preserves the coloring \(f\) on \(L(T)\). By the definition of the partial coloring \(g_\varepsilon\) on the path \(uvw_2\), both \(uv,vw_2\) are bichromatic. Therefore, there are at most two newly created monochromatic edges, one edge joining \(w_2\) and a leaf, and possibly \(vw_1\). Thus, the set of monochromatic edges remains a matching in \(T\), since the retained monochromatic edges form a matching in \(T^{**}-\{w_1\}\), and the two possible new edges have disjoint endvertices, none incident with a retained monochromatic edge. So, conclusions (1) and (2) are verified.

Now we are ready to verify conclusion (3). The vertices of \(T\) that are not in \(T^{**}\) are \(v,w_2\) and the two leaf neighbors of \(w_2\). The vertices \(v,w_2\) have opposite colors, and the two leaf neighbors of \(w_2\) also have opposite colors. Thus, the difference between \(|B|\) and \(|W|\) under \(f'_\varepsilon\) remains the same as that under \(h'_\varepsilon\). That is, for each \(\varepsilon\in\{+1,-1\}\),
\[
D_f=D_h,\qquad D_{f'_\varepsilon}=D_{h'_\varepsilon}.
\]
By induction, the bounds in \eqref{eq:tree-bounds} hold for \(h\) and \(h'_\varepsilon\), that is,
\[
\frac{2D_h-2}{3}\le D_{h'_{+1}}\le\frac{2D_h+6}{3},
\qquad
\frac{2D_h-6}{3}\le D_{h'_{-1}}\le\frac{2D_h+2}{3}.
\]
Combining these inequalities with the equalities above verifies \eqref{eq:tree-bounds}. That is,
\[
\frac{2D_f-2}{3}\le D_{f'_{+1}}\le\frac{2D_f+6}{3},
\qquad
\frac{2D_f-6}{3}\le D_{f'_{-1}}\le\frac{2D_f+2}{3}.
\]
\end{proof}

\section{Switching along bichromatic paths}

For a vertex \(2\)-coloring of a cubic graph, a vertex is of \emph{type~\(0\)}, \emph{type~\(1\)}, \emph{type~\(2\)}, or \emph{type~\(3\)} if it has exactly \(0,1,2\), or \(3\) neighbors of its own color, respectively.

Let \(A_f\) be the subgraph of \(G\) induced by the set of bichromatic edges. Every path in \(A_f\) is alternately colored. Vertices of types \(0,1,2\) are of degrees \(3,2,1\) in \(A_f\), respectively. Notice that vertices of types \(0\) and \(2\) are of odd degree in \(A_f\).

We first show the existence of a path from a type~\(2\) vertex to another vertex of odd degree in \(A_f\).

\begin{obs}\label{obs:odd-path}
Let \(f\) be a vertex \(2\)-coloring of a cubic graph \(G\), and let \(v_1\) be a type~\(2\) vertex. Let \(Q_1\) be the component of \(A_f\) containing \(v_1\). Then there is a shortest path \(P=v_1\cdots v_t\) from \(v_1\) to the set of all other vertices in \(V(Q_1)\) with odd degree in \(A_f\). Furthermore, every internal vertex of \(P\) is of type~\(1\), and \(v_t\) is of type~\(0\) or type~\(2\).
\end{obs}
\begin{proof}
Since \(v_1\) is of type~\(2\), it is of odd degree in \(A_f\). The component \(Q_1\) contains an even number of vertices with odd degree in \(A_f\), so there is another such vertex. This proves that such a shortest path exists. By the choice of \(P\), every internal vertex is of even degree in \(A_f\) and is incident with two edges of \(P\). Since \(G\) is cubic, each internal vertex has exactly one monochromatic neighbor and is therefore of type~\(1\). The terminal vertex is of odd degree in \(A_f\), so it is of type~\(0\) or type~\(2\).
\end{proof}

\emph{Switching} a vertex set \(X\) interchanges the colors of its vertices and leaves all other colors unchanged. Switching a path \(P\) means switching \(V(P)\).

\begin{obs}\label{obs:switched-path}
Let \(f\) be a vertex \(2\)-coloring of a cubic graph \(G\), and let \(P'=v_1\cdots v_r\) be a path in \(A_f\). Let \(f'\) be the coloring obtained by switching \(V(P')\) in \(f\). Then every vertex in \(V(P')\setminus\{v_1,v_r\}\) is of type~\(0\) or type~\(1\) under \(f'\).
\end{obs}
\begin{proof}
Both ends of every edge in \(E(P')\) are switched, so these edges remain bichromatic. Each vertex in \(V(P')\setminus\{v_1,v_r\}\) therefore is incident with at most one monochromatic edge under \(f'\).
\end{proof}

We apply a switching operation to a subpath of the path \(P=v_1\cdots v_t\) described in Observation~\ref{obs:odd-path}. Without loss of generality, we assume that the first vertex \(v_1\) of \(P\) is colored black.

\begin{figure}[H]
\centering
\includegraphics[width=0.96\linewidth]{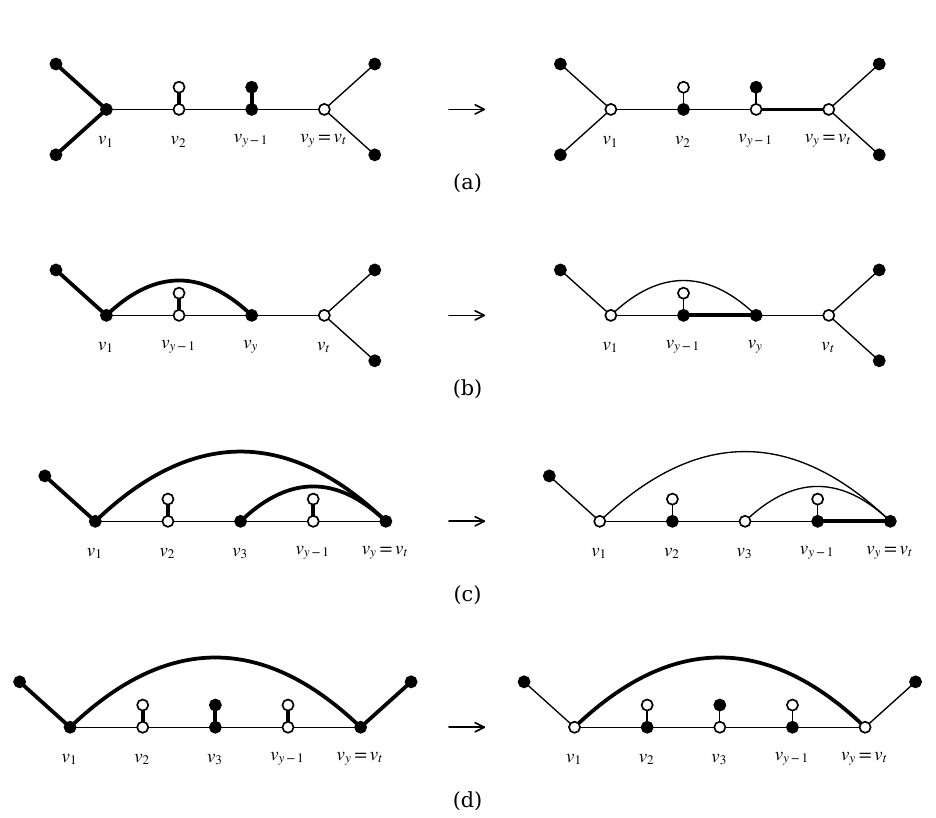}
\caption{Typical switching operations in Lemma~\ref{result:4.1}, before (left) and after (right). In (a)--(c), \(P'=Q\), and \(v_y\) is originally of type~\(0\), type~\(1\), and type~\(2\), respectively. In (c), both monochromatic edges at \(v_y\) join it to \(V(Q)\). In (d), \(P'=P\), and the chord \(v_1v_t\) remains monochromatic. Bold edges are monochromatic.}\label{fig:2}
\end{figure}

\begin{lem}\label{result:4.1}
Let \(f\) be a vertex \(2\)-coloring of a cubic graph \(G\) with no type~\(3\) vertices. Choose \(P=v_1\cdots v_t\) as in Observation~\ref{obs:odd-path}, with \(v_1\) black. Then there exists a subpath \(P'=v_1\cdots v_r\) of \(P\), where \(1\le r\le t\), such that switching \(V(P')\) yields a new coloring with the following conclusions:
\begin{enumerate}[label=(\arabic*)]
\item every original type~\(2\) vertex of \(N_G[V(P')]\) becomes type~\(0\) or type~\(1\);
\item no new type~\(2\) vertex and no type~\(3\) vertex are created; and
\item the number of black vertices does not increase and decreases by at most \(1\).
\end{enumerate}
\end{lem}
\begin{proof}
By Observation~\ref{obs:odd-path}, every internal vertex of \(P\) is of type~\(1\), and \(v_t\) is of type~\(0\) or type~\(2\).
Recall that a chord of \(P\) in \(G\) is an edge of \(G\) joining two nonconsecutive vertices in \(V(P)\). Now, we construct the subpath \(P'\) step by step as follows.

\textbf{Step 1.} If \(P\) has a monochromatic chord \(v_iv_j\) in \(G\) with \(1\le i<j\le t\), let \(y\) be the smallest possible value of \(j\). If \(P\) has no monochromatic chord, set \(y=t\).

Let \(Q=v_1\cdots v_{y-1}\) be the corresponding subpath of \(P\).

\textbf{Step 2.} If \(v_y\) is of type~\(2\) and at most one monochromatic edge joins \(v_y\) to \(V(Q)\), set \(P'=v_1\cdots v_y=Qv_y\). Otherwise, set \(P'=v_1\cdots v_{y-1}=Q\).

Before switching \(V(P')\), we present two properties of \(Q\).

\textbf{Claim 1.} Every edge of \(E(G)\setminus E(P)\) incident with \(V(Q)\) is monochromatic.

Vertex \(v_1\) is of type~\(2\), so its only bichromatic edge belongs to \(E(P)\). By Observation~\ref{obs:odd-path}, every vertex in \(V(Q)\setminus\{v_1\}\) is of type~\(1\) and is incident with two bichromatic edges of \(E(P)\). Thus, its remaining incident edge is monochromatic. This proves Claim 1.

\textbf{Claim 2.} \(Q\) is a chordless path.

For otherwise, \(Q\) would have a chord \(v_iv_j\) with \(i<j<y\). By Claim 1, this chord would be monochromatic, contradicting the choice of \(y\) in Step 1.

Now, we apply the switching operation to \(V(P')\).

The following claim is straightforward.

\textbf{Claim 3.} Every edge in \(E(P')\) remains bichromatic. An edge changes its monochromatic or bichromatic status if and only if it has exactly one endvertex in \(V(P')\).

We now verify the three conclusions.

For conclusions (1) and (2), we consider the two possible choices of \(P'\) (see Step 2).

\medskip
\textbf{Case 1.} \(V(P')=V(Q)\). The edge \(v_{y-1}v_y\) is the only edge in \(E(P)\) with exactly one endvertex in \(V(Q)\), and it becomes monochromatic by Claim 3. By Claims 1 and 2, every edge of \(E(G)\setminus E(P)\) incident with \(V(Q)\) is originally monochromatic and has exactly one endvertex in \(V(Q)\). All these edges become bichromatic by Claim 3. Thus \(v_{y-1}v_y\) is the only new monochromatic edge.

We now determine the types of the vertices in \(V(Q)\) after the switching operation. If \(y>2\), every edge incident with \(v_1\) is bichromatic, so \(v_1\) becomes type~\(0\). The vertex \(v_{y-1}\) becomes type~\(1\), since its only monochromatic edge is \(v_{y-1}v_y\). If \(y=2\), then \(Q\) consists of the single vertex \(v_1=v_{y-1}\), which becomes type~\(1\). The internal vertices of \(Q\) are of type~\(0\) or type~\(1\) by Observation~\ref{obs:switched-path}.

We next consider the vertices in \(V(G)\setminus(V(Q)\cup\{v_y\})\). Their monochromatic degrees cannot increase. In particular, every original type~\(2\) vertex in \(N_G[V(Q)]\setminus(V(Q)\cup\{v_y\})\) loses at least one monochromatic edge and gains none. Each such vertex becomes type~\(0\) or type~\(1\).

Only \(v_y\) remains to be checked. If it is originally of type~\(0\) (see Figure~\ref{fig:2}(a)), then it becomes type~\(1\), since \(v_{y-1}v_y\) is its only monochromatic edge after the switching operation. If \(v_y\) is originally of type~\(1\) (see Figure~\ref{fig:2}(b)), then \(v_y\) is of degree two in \(A_f\). By the choice of \(P\) in Observation~\ref{obs:odd-path}, we have \(y<t\). By Step 1, there is a monochromatic chord joining \(v_y\) to \(V(Q)\). This chord is the only monochromatic edge incident with \(v_y\) under \(f\). The switching operation makes this chord bichromatic and \(v_{y-1}v_y\) monochromatic, so \(v_y\) remains type~\(1\). If \(v_y\) is originally of type~\(2\) (see Figure~\ref{fig:2}(c)), then \(v_y\) is of degree one in \(A_f\). By the choice of \(P\) in Observation~\ref{obs:odd-path}, we have \(y=t\). Step 2 gives \(P'=Q\) only when both monochromatic edges incident with \(v_y\) join it to \(V(Q)\). Both become bichromatic, and \(v_{y-1}v_y\) becomes the only monochromatic edge incident with \(v_y\). Thus \(v_y\) is of type~\(1\) after the switching operation in all three cases. This proves (1) and (2) when \(V(P')=V(Q)\).

\medskip
\textbf{Case 2.} \(V(P')=V(Q)\cup\{v_y\}\). By Step 2, \(v_y\) is of type~\(2\). Since \(v_1,v_t\) are the only possible type~\(2\) vertices in \(V(P)\), we have \(y=t\) and \(P'=P\). Thus, by Claim 1, \(P\) is an \(A_f\)-component, and every edge with exactly one endvertex in \(V(P)\) is monochromatic under \(f\). By Claim 2, every monochromatic chord of \(P\) is incident with \(v_t\). Step 2 allows at most one such chord. Figure~\ref{fig:2}(d) shows the switching operation when the chord is \(v_1v_t\).

Switching \(V(P)\) keeps the edges in \(E(P)\) bichromatic and the possible chord monochromatic. Every edge with exactly one endvertex in \(V(P)\) becomes bichromatic by Claim 3. The vertex \(v_1\) becomes type~\(1\) if \(v_1v_t\) is a chord of \(P\), and type~\(0\) otherwise. Since every chord of \(P\) is incident with \(v_t\), the vertex \(v_t\) becomes type~\(1\) if a chord exists, and type~\(0\) otherwise. The internal vertices are of type~\(0\) or type~\(1\) by Observation~\ref{obs:switched-path}.

Every vertex in \(N_G[V(P)]\setminus V(P)\) loses at least one monochromatic edge and gains none. Thus no new type~\(2\) or type~\(3\) vertex is created, and every original type~\(2\) vertex in this set becomes type~\(0\) or type~\(1\). Vertices in \(V(G)\setminus N_G[V(P)]\) retain their types. This proves (1) and (2).

\medskip
Finally, we check conclusion (3) for the resulting coloring. Since the colors alternate along \(P'\) and \(v_1\) is black, \(V(P')\) contains equally many vertices of each color or one more black vertex. Switching \(V(P')\) therefore leaves the number of black vertices unchanged or decreases it by \(1\).
\end{proof}

In the following lemma, we apply Lemma~\ref{result:4.1} repeatedly to reduce the number of type~\(2\) vertices.

\begin{lem}\label{result:4.2}
Let \(f\) be a vertex \(2\)-coloring of a cubic graph \(G\) with no type~\(3\) vertices, and suppose that every type~\(2\) vertex is black. Let \(H\) be the union of the monochromatic components of order at least three, and \(R\subseteq V(H)\) such that every component of \(H-R\) is of order at most two. Then there is a new coloring \(f'\) such that every vertex is of type~\(0\) or type~\(1\), while the number of black vertices does not increase and decreases by at most \(|R|\).
\end{lem}
\begin{proof}
We proceed by induction on \(|R|\). If \(|R|=0\), every monochromatic component is of order at most two, so \(f\) itself satisfies the conclusion. Suppose that \(|R|\ge1\) and that the statement holds for all smaller sets satisfying the hypotheses.

Let \(r\in R\). Since \(r\in V(H)\), it is of type~\(1\) or type~\(2\). If \(r\) is of type~\(1\), its unique monochromatic neighbor is of type~\(2\); let this neighbor be \(v_1\). If \(r\) is of type~\(2\), let \(v_1=r\). In either case, \(v_1\) is black.

Observation~\ref{obs:odd-path} gives a path starting at \(v_1\) to which Lemma~\ref{result:4.1} applies. Let \(P'\) be the subpath given by the lemma, and let \(g\) be the coloring obtained by switching \(V(P')\). Every original type~\(2\) vertex in \(N_G[V(P')]\) becomes type~\(0\) or type~\(1\), and no new type~\(2\) or type~\(3\) vertex is created. Any remaining type~\(2\) vertices belong to \(V(G)\setminus N_G[V(P')]\) and remain black. The number of black vertices does not increase and decreases by at most \(1\).

Let \(R'=R\setminus\{r\}\), and let \(H'\) be the union of the monochromatic components of order at least three under \(g\). To apply the induction hypothesis, we still need to verify that \emph{every component of \(H'-R'\) is of order at most two}. To verify this statement, it suffices to prove that every path \(xvz\) of order three in \(H'\) satisfies
\begin{equation}\label{eq:remaining-R}
\{x,v,z\}\cap R'\ne\varnothing.
\end{equation}
Consider any such path. Its middle vertex \(v\) is of type~\(2\) under \(g\). By conclusions (1) and (2) of Lemma~\ref{result:4.1}, \(v\) was of type~\(2\) under \(f\) and lies outside \(N_G[V(P')]\). Hence \(\{x,v,z\}\cap V(P')=\varnothing\). Thus the colors of \(x,v,z\) remain unchanged, and \(xvz\) was also a monochromatic path in \(H\). Therefore, \(\{x,v,z\}\cap R\ne\varnothing\). To verify \eqref{eq:remaining-R}, it remains to show that
\begin{equation}\label{eq:avoid-r}
\{x,v,z\}\cap\{r\}=\varnothing.
\end{equation}
Suppose otherwise. If \(r\) was of type~\(1\) under \(f\), then it is of type~\(0\) or type~\(1\) under \(g\) by conclusion (2) of Lemma~\ref{result:4.1}. Since \(r\in\{x,v,z\}\), it must be of type~\(1\) under \(g\) and an endvertex of \(xvz\). The path \(xvz\) is also monochromatic under \(f\), where the unique monochromatic neighbor of \(r\) is \(v_1\). Hence \(v=v_1\), contradicting the fact that \(v_1\) is no longer of type~\(2\) under \(g\). Thus \eqref{eq:avoid-r} holds in this case.

If \(r\) was of type~\(2\) under \(f\), then \(v_1=r\), so \(r\in V(P')\). The path \(xvz\) is also monochromatic under \(f\). Since \(r\in\{x,v,z\}\), the middle vertex \(v\) is either \(r\) itself or a neighbor of \(r\). Thus \(v\in N_G[r]\subseteq N_G[V(P')]\), contradicting \(v\notin N_G[V(P')]\). Thus \eqref{eq:avoid-r} holds in this case.

Since \(H'-(R'\cap V(H'))=H'-R'\), the set \(R'\cap V(H')\) satisfies the required condition. We can now apply the induction hypothesis to \(g\) and this set. This gives a coloring \(f'\) in which every vertex is of type~\(0\) or type~\(1\). The number of black vertices does not increase, and the further decrease is at most \(|R'\cap V(H')|\le |R'|=|R|-1\). Including the first switching operation, the total decrease is at most \(1+(|R|-1)=|R|\).
\end{proof}

\section{Proof of the main theorem}

We prove Theorem~\ref{result:1.1} in this section. The following corollary of Lemma~\ref{result:4.2} is a main ingredient in the proof.

\begin{cor}\label{result:5.1}
Let \(k\in\{1,2\}\), and let \(f:V(G)\to\{B,W\}\) be a vertex \(2\)-coloring of a cubic graph \(G\), with \(|B|-|W|=2k-2\). Suppose that the monochromatic edge set under \(f\) can be partitioned into two matchings \(M_1\) and \(M_2\), where \(M_2\) consists of \(k\) black monochromatic edges. Then \(G\) admits a \(2\)-quasi-bisection.
\end{cor}

\begin{proof}
To apply Lemma~\ref{result:4.2}, we are required to specify \(H\) and a suitable set \(R\). Let \(H\) be the union of the monochromatic components of order at least three. Since \(M_1\) and \(M_2\) are matchings, no vertex is of type~\(3\), and the edges in each component of \(H\) alternate between \(M_1\) and \(M_2\). Every type~\(2\) vertex is incident with an edge of \(M_2\) and is therefore black. Let \(R\) consist of one chosen endvertex of each edge in \(M_2\cap E(H)\). Then \(R\subseteq V(H)\) and \(|R|\le |M_2|=k\). Every edge of \(H-R\) belongs to \(M_1\), so every component of \(H-R\) is of order at most two.

By applying Lemma~\ref{result:4.2}, we obtain a new coloring with only type~\(0\) and type~\(1\) vertices, in which the number of black vertices does not increase and decreases by at most \(k\). Since the initial color-class difference is \(2k-2\) and \(k\in\{1,2\}\), the resulting color classes \(B',W'\) satisfy
\[
-2\le |B'|-|W'|\le 2k-2\le2.
\]
\end{proof}

We now prove Theorem~\ref{result:1.1} by applying Corollary~\ref{result:5.1} and Lemma~\ref{result:3.2}.

\begin{proof}[Proof of Theorem~\ref{result:1.1}]
Let \(T\) be the tree in the statement and \(J\) the subgraph induced by \(E(G)\setminus E(T)\). Since \(J\) is even, Observation~\ref{obs:even-complement-tree} implies that \(T\) is a spanning cubic tree of \(G\). Each leaf of \(T\) is incident with two edges of \(J\), and each internal vertex is incident with none. Hence every vertex of \(J\) has degree two. Thus \(J\) is a disjoint union of circuits and \(V(J)=L(T)\). Let \(k\) be the number of odd circuit components of \(J\). If \(k=0\), then \(J\) is bipartite, and the result follows from DeVos and Nurse \cite{ref4}.

We may therefore assume that \(k\in\{1,2\}\). Let \(C_1,\ldots,C_k\) be the odd circuit components of \(J\). First, we construct a coloring \(f\) of \(V(J)\) as follows. Color each odd circuit \(C_i\) so that it has exactly one black monochromatic edge \(e_i\) and all other edges are bichromatic. Color every even circuit alternately black and white. Then \(D_f=k\).

Since \(V(J)=L(T)\), we can apply Lemma~\ref{result:3.2} to extend the partial coloring \(f\) to a new coloring \(f':V(G)\to\{B,W\}\), taking
\[
f'=
\begin{cases}
f'_{-1}, & k=1,\\
f'_{+1}, & k=2.
\end{cases}
\]
To obtain a \(2\)-quasi-bisection from \(f'\), we apply Corollary~\ref{result:5.1}. We first identify the matchings \(M_1,M_2\) required by the corollary. The monochromatic tree edges under \(f'\) form a matching \(M_1\). Since the leaf colors are preserved, the monochromatic edges of \(J\) under \(f'\) are exactly \(M_2=\{e_1,\ldots,e_k\}\). These edges lie in distinct circuit components of \(J\), so \(M_2\) is a matching of \(k\) black monochromatic edges. As \(E(T)\) and \(E(J)\) partition \(E(G)\), the monochromatic edge set under \(f'\) is the disjoint union of the matchings \(M_1\) and \(M_2\).

Since \(G\) is cubic, \(|B|+|W|=|V(T)|=|V(G)|\) is even. Therefore, so is \(|B|-|W|\). Thus, by equation \eqref{eq:tree-bounds}, we have
\[
|B|-|W|=
\begin{cases}
0, & k=1,\\
2, & k=2,
\end{cases}
\]
that is, \(|B|-|W|=2k-2\) for \(k\in\{1,2\}\). Corollary~\ref{result:5.1} can therefore be applied to complete the proof.
\end{proof}

\section*{Acknowledgments}

This work was partially supported by the National Natural Science Foundation of China (Nos. 12471321 and 12331013).

\begingroup\small\raggedright

\section*{Declaration of interests and the data availability}

We declare that we have no known competing financial interests or personal relationships that
could have appeared to influence the work reported in this paper.\\

No data were generated or analysed in this manuscript.

\endgroup

\begin{thebibliography}{99}
\bibitem{ref1} M. Abreu, J. Goedgebeur, D. Labbate, and G. Mazzuoccolo, A note on \(2\)-bisections of claw-free cubic graphs, \emph{Discrete Applied Mathematics} \textbf{244} (2018), 214--217. \href{https://doi.org/10.1016/j.dam.2018.03.016}{doi:10.1016/j.dam.2018.03.016}.
\bibitem{ref2} M. Abreu, J. Goedgebeur, D. Labbate, and G. Mazzuoccolo, Colourings of cubic graphs inducing isomorphic monochromatic subgraphs, \emph{Journal of Graph Theory} \textbf{92} (2019), 415--444. \href{https://doi.org/10.1002/jgt.22462}{doi:10.1002/jgt.22462}.
\bibitem{ref3} A. Ban and N. Linial, Internal partitions of regular graphs, \emph{Journal of Graph Theory} \textbf{83} (2016), 5--18. \href{https://doi.org/10.1002/jgt.21909}{doi:10.1002/jgt.21909}.
\bibitem{ref4} M. DeVos and K. Nurse, \emph{On the Ban--Linial Conjecture}, \href{https://arxiv.org/abs/2512.18913}{arXiv:2512.18913} (2025).
\bibitem{ref5} L. Esperet, G. Mazzuoccolo, and M. Tarsi, Flows and bisections in cubic graphs, \emph{Journal of Graph Theory} \textbf{86} (2017), 149--158. \href{https://doi.org/10.1002/jgt.22118}{doi:10.1002/jgt.22118}.
\bibitem{ref6} X. Li, R.-X. Hao, Y. Luo, and S. Liu, Ban--Linial's Conjecture and Halin snarks, \emph{Journal of Graph Theory} \textbf{113} (2026), 527--537. \href{https://doi.org/10.1002/jgt.70080}{doi:10.1002/jgt.70080}.
\bibitem{ref8} J. P. Zerafa, Ban--Linial's Conjecture and treelike snarks, \emph{Journal of Graph Theory} \textbf{101} (2022), 511--520. \href{https://doi.org/10.1002/jgt.22837}{doi:10.1002/jgt.22837}.
\end{thebibliography}
\end{document}